\documentclass[12pt]{amsart}
\usepackage{amsfonts,amsmath,amsthm,oldgerm,amssymb,amscd}

\usepackage[all,cmtip]{xy}
\newcommand{\Q}{{\mathbb Q}}
\newcommand{\R}{{\mathbb R}}
\newcommand{\Z}{{\mathbb Z}}
\newcommand{\C}{{\mathbb C}}

\newcommand{\A}{{\mathbf A}}
\newcommand{\G}{{\mathbf G}}

\renewcommand{\H}{{\mathbb H}}

\newtheorem{theorem}{Theorem}
\newtheorem{corollary}{Corollary}
\newtheorem{proposition}[theorem]{Proposition}

\theoremstyle{definition}
\newtheorem{conjecture}{Conjecture}
\newtheorem*{example}{Example}
\newtheorem{definition}{Definition}
\numberwithin{definition}{section}
\newtheorem{question}{Question}
\numberwithin{equation}{section}

\theoremstyle{remark}
\newtheorem{remark}{Remark}
\newtheorem*{acknowledgement}{Acknowledgement}

\begin{document}
\title[On spectrum and arithmetic]{Some questions on spectrum and
  arithmetic\\of locally symmetric spaces}

\author[C. S. Rajan]{C. S. Rajan}

\noindent\thanks{Journal Ref: Advanced Studies in Pure Mathematics 58, 2010\\
Algebraic and Arithmetic Structures of Moduli Spaces (Sapporo
2007)\\pp 137-157}

\address{
Tata Institute of Fundamental Research\\
Homi Bhabha Road\\ 
Mumbai 400 005\\
India}
\email{rajan@math.tifr.res.in}

\subjclass[2000]{}


\begin{abstract}
We consider the question that the  spectrum and arithmetic of
locally symmetric spaces defined by congruent arithmetical lattices
should mutually determine each other. We frame these questions in the
context of automorphic representations.  
\end{abstract}

\maketitle

\section{Introduction}  The classical spectral theories of light and
sound can be mathematically interpreted and generalized in possibly
two different ways: one, in the context of Riemannian manifolds and
the spectrum of the Laplacian acting on a suitably defined space of
smooth functions on the manifold. More generally, one can consider the
spectrum of elliptic, self-adjoint, natural differential operators acting on 
subspaces defined by appropriate boundary conditions of the space of
smooth sections of natural vector bundles on the manifold, for
example, the deRham Laplacian acting on the space of smooth forms. 

The concept of spectrum can also be generalized in the context of
continuous representation theory of topological groups. In what
follows, we restrict our attention to the class of locally symmetric
Riemannian manifolds.
 These are spaces of the form $\Gamma\backslash G/K$, where $G$
is a reductive Lie group, $K$ is a maximal compact subgroup of $G$,
and $\Gamma$ is a lattice in $G$. The universal cover $G/K$ carries a
natural $G$-invariant metric, the Bergman metric, and we can relate
the representation theory of $G$ to the spectrum of Laplace type
operators on $\Gamma\backslash G/K$.

It is well known that the spectra  of Laplace type operators
associated with locally symmetric Riemannian manifolds or the
representation theory of $G$ acting on the space of functions on
$\Gamma\backslash G$, share many similarities with the arithmetic of
such spaces. The beginnings of this analogy can be attributed to Maass,
who introduced and studied the arithmetic properties   of
eigenfunctions  of the Laplace operator on the upper half plane which
are invariant under a congruence subgroup of $SL_2(\Z)$ (called Maass
forms) \cite{Ge}.  
It was shown that Hecke operators can be defined on the space
of such eigenfunctions of the Laplacian, and that the $L$-functions
associated to Maass-Hecke eigenforms have many similar properties to the
$L$-functions associated to holomorphic Hecke eigenforms. 

The theory
of Maass forms also provides a link between the representation theory
of $G=PSL_2(\R)$ acting on the space of square integrable functions
$L^2(\Gamma\backslash G)$, and the spectrum of the hyperbolic
Laplacian acting on $\Gamma \backslash \H$, where  $\H$ is hyperbolic
upper half plane. For general groups, Langlands theory of automorphic
representations provides the appropriate generalization to relate not
only the spectrum of invariant differential operators to
representation theory, but also to the arithmetic that is provided by
this theory. 

It was Selberg who brought the analogy between the spectrum
and the arithmetic to the
fore. Selberg considered primitive closed
geodesics on a hyperbolic surface of finite area, as analogous to
rational primes, and established an analogue of the prime number
theorem for primitive closed geodesics \cite{He}. The concept of the length
spectrum can be introduced, and the Selberg trace formula  establishes
a relationship between the spectrum of the Laplacian acting on
funtions on a compact hyperbolic surface and the length spectrum of
primitive closed geodesics on it. Although the relationship with 
the length spectrum is quite important, we will not pursue it 
out here. 

Another analogy observed by  Selberg, is the
conjecture that
\[ \lambda_1(\Gamma \backslash \H)\geq 1/4,\]
where $\lambda_1(\Gamma \backslash \H)$ is the least non-zero
eigenvalue of the Laplacian acting on the space of smooth functions on
$\Gamma \backslash \H$. Here $\Gamma$ is a congruent arithemetical
lattice contained inside $SL_2(\Z)$. 
Selberg's conjecture can be considered as the archimedean
analogue of Ramanujan's conjecture on the size of the Hecke
eigenvalues of holomorphic newforms \cite{Sa}.  

In a similar vein, there are many common features between the
representation theories of real and $p$-adic Lie groups. For example,
the Harishchandra homomorphism identifies the center $Z(\frak g)$ of
the universal enveloping algebra of a complex semisimple Lie algebra $\frak g$
with the ring of Weyl group invariants of the symmetric algebra 
of a maximal torus in $\frak g$ \cite{Kn}. The Casimir elements belong to
$Z(\frak g)$, and these project to invariant differential operators
on the associated symmetric space.   The  non-archimedean counterpart
of the Harishchandra isomorphism is the Satake isomorphism:  if $G$ is
a split reductive $p$-adic Lie group,   $K$ a hyperspecial maximal
compact subgroup of $G$, $T$ a maximal split torus of $G$, the Satake
isomorphism identifies the Hecke algebra $H(G,K)$ of  $K$-biinvariant
functions on $G$,  to the Weyl group invariants of the algebra of
unramified characters of $T$ (see Cartier's article in \cite{Cor}).
 Thus the Hecke algebra 
is analogous to the algebra of invariant
differential operators of a complex semisimple Lie algebra.

The inverse spectral problem is to know the properties of the
Riemannian manifold that are determined by the spectrum. For instance,
is the isometry class of the space determined by the spectrum? In the
context of planar domains, this question was posed in a colloquiual
way by M. Kac as ``Can one hear the shape of the drum?''.
Milnor and later Vigneras used arithmetical
methods to  construct pairs of isospectral but non-isometric
manifolds. A fundamental construction of pairs of isospectral
manifolds was given by Sunada. Sunada's method is completely analogous
to Gassmann's method  of exhibiting 
non-conjugate number fields having the same Dedekind zeta function.  
These examples indicate the existence of an arithmetical flavor to the
inverse spectral problem. 

Based on such examples, a well known heuristic is that the Laplacian
can be considered as the Frobenius/Hecke  operator at infinity, for
the class of locally symmetric spaces defined by congruent
arithmetical lattices. Such spaces have a natural Riemannian structure
arising from the invariant metric on the universal cover.   When the
universal cover is hermitian symmetric, a natural arithmetical
structure is provided by  the theory of canonical models due to
Shimura and Deligne. In many examples, these spaces are the moduli
spaces of abelian varieties endowed with extra structure. 
In the general case, the  Langlands theory of
automorphic representations can be considered as an automorphic aspect
of the arithmetic of such spaces.  

We expound on the  theme that the
Frobenius at infinity is the Laplacian,  and our  basic expectation is
that the spectrum and arithmetic for such spaces should mutually
determine each other.

\begin{acknowledgement} My sincere thanks to Dipendra Prasad for many 
discussions on this topic.  I thank T. Chinburg,
L. Clozel, C. Deninger, S. Ganguly, G, Harder, S. Kudla, M. S. Narasimhan, 
G. Prasad, N. Ramachandran, 
A. W. Reid, P. Sarnak and A. Venkatesh for many useful 
discussions, remarks, etc. 

During the course of this work, I availed of the hospitality of MPIM
at Bonn during the year 2005-06 and April of 2009, the IAS at
Princeton for the first half of 2007, and the University of M\:{u}nster
during May-June of 2008. I thank all these Institutes for their
support and wonderful hospitality.
\end{acknowledgement}
   
\section{Archimedean analogue of Tate's conjecture}
Let $M$ be a compact Riemann  surface, uniformized by the upper half
plane.  The  Jacobian $J(M)$ of $M$, the `motive'
associated to the Riemann surface $M$, can be considered as the
`arithmetic' associated to $M$. 

We recall the theorem of Faltings proving Tate's
conjecture \cite{F}: if the eigenvalues of the Frobenius elements
acting on the $l$-adic cohomology groups of two smooth, projective
curves defined over a number field coincide, then the Jacobians of the
curves are isogenous. The analogy of the Laplacian with the Frobenius
motivates us to ask a naive archimedean analogue of Tate's conjecture:

\begin{question}\label{qn:naiveaat}
Suppose $M_1, ~M_2$ are compact Riemann surfaces of genus at least
two, that are isospectral
with respect to the hyperbolic Laplacian acting on the space of smooth
functions. Are the Jacobians $J(M_1)$ and  $J(M_2)$ isogenous? 
\end{question}
This question turns out to have a false answer (see Corollary
\ref{cor:naiveatfalse})
in general, but can be verified for the examples constructed using the
Gassmann-Sunada method (\cite{PR}). 

\begin{remark} A similar question was raised by M. S. Narasimhan (to
  T. Sunada) but requiring the stronger conclusion that the Jacobians
  be isomorphic. This seems unrealistic to expect even for the
  examples constructed using the Sunada method. From the arithmetic
  viewpoint it is natural to consider the isogeny class.  

Narasimhan also raised the quesion of considering isospectrality with
respect to other natural metrics defined on hyperbolic compact Riemann
surfaces, for example the pullback of the Fubini-Study metric with
respect to the canonical map to projective space, or the pullback of
the translation invariant metric with respect to embedding the curve
in it's Jacobian, etc. But it is the hyperbolic metric that is natural
in the context of automorphic forms (see 
Theorem \ref{thm:convsunadpsl2}), and thus
seems natural to consider if we want to relate the spectrum 
and the arithmetic of such spaces. 
\end{remark}

\subsection{Flat tori} 
By the theorem of Faltings, we see that two abelian varieties defined
over a number field have isomorphic $l$-adic representations if and
only if they are isogenous. 

In the context of flat tori, define two flat tori to be {\em arithmetically
equivalent} if they are isogenous. We define two flat tori  to be
{\em spectrally commensurable} if upto rescaling of metrics the
spectra are mutually contained inside each other.
The expectation that the spectrum determines the arithmetic
  translates to the following refinement of a conjecture of Kitaoka \cite{K}:
\begin{conjecture} Suppose $L_1, ~L_2$ are lattices in Euclidean
  $n$-space $\R^n$, such that the corresponding  flat tori $\R^n/L_1,
  ~\R^n/L_2$ are spectrally commensurable. 
Then upto an isometry of $\R^n$, the lattices $L_1$ and $L_2$ are
commensurable. 
\end{conjecture}
(To put it more colloquiually, the thetas and the zetas mutually
determine each other.) 

The  known examples (see \cite{Con} for some of them),
 for example Milnor's example of the isospectral
pair of flat tori given by the lattices $(E_{16}, E_8\oplus E_8)$, can
be seen to be commenurable. Indeed, 
when the quadratic form associated to the lattices are integral
valued, a brief hint of the proof of the above conjecture with the
stronger assumption that the flat tori are isospectral,  was given  by
Kitaoka \cite{K}, and a detailed  proof of a more general result was
given by Bayer and Nart \cite{BN}.

\section{Arithmetic Archimedean Analogue of Tate's conjecture}
\label{section:aaat}
In the previous section, we had taken as our model for the arithmetic
of a compact Riemann surface, the Jacobian of the surface. Here, we 
focus our attention on the example of projective Shimura curves: they
have natural canonical models defined over a number field, and so we
can consider the usual arithmetic of such spaces \cite{Sh2}.

Let $F$ be a number field and let $D$ be  an indefinite,  quaternion
division algebra  over $F$. Let 
\[   G=SL_1(D), \quad \mbox{and}\quad \tilde{G}=GL_1(D).\]  
Let $K_{\infty}$ be a maximal compact subgroup of
$G_{\infty}=G(F\otimes \R)$ and let $M=G_{\infty}/K_{\infty}$ be the
non-compact symmetric space associated to $G$. If $D$ is unramified at
$r$ real places and $2r_2$ is the number of complex embeddings of $F$,
then, 
\[ M\simeq \H_2^r\times \H_3^{r_2},\]
where $\H_2$ (resp. $\H_3$) is the hyperbolic $2$-space
(resp. $3$-space). 
For a number field $F$,
let $\A_F$ denote the
adele ring of $F$, and $\A_{F,f}$ the subring of finite adeles. Let $K$ be a
compact open subgroup of $G(\A_{F,f})$. Let  $\Gamma_K$ be  the co-compact
lattice  in $G_{\infty}$ corresponding to $K$, given by the projection
to $G_{\infty}$ of the group  $G_{\infty}K\cap G(F)$.   If $K$ is
sufficiently small, then  $\Gamma_K$ will be a torsion-free lattice
acting freely and properly discontinuously on $M$.  The quotient space
$M_K=\Gamma_K\backslash M$ will then be  a compact, locally symmetric
space.

Suppose now  $F$ is totally real. We recall the theory of canonical
models due to Shimura \cite{Sh1}. Let $\tau_1, \cdots, \tau_r$
be the real embeddings corresponding to archimedean places of $F$ at
which $D$ splits. The reflex field $F'$ of $(F, \tau_1,
\cdots, \tau_r)$ is the field generated by the sums
\[ \sum_{i=1}^r\tau_i(x), ~x\in F.\] 
There exists a $\Q$-rational homomorphism 
\[ \lambda:~ R_{F'/\Q}\G_m\to R_{F/\Q}\G_m, \]
where for a number field $K$, $R_{K/\Q}$ denotes the Weil restriction
of scalars. On the $\Q$-rational points, $\lambda:F'^*\to F^*$
satisfies the following relation:
\[ N_{F/\Q}(\lambda(x))=N_{F'/\Q}(x)^r, \quad x\in F'^*.\]
When $r=1$, the reflex field can be taken to be $F$ itself. 

Let $\tilde{G}_{\infty +}$ (resp. $F_{\infty +}^*$)  be the identity component of 
 $\tilde{G}_{\infty }$ (resp. $F_{\infty }^*$). Let
\[ \tilde{\mathcal G}_+=\{x\in \tilde{G}_{\infty +}\tilde{G}(\A_{F,f})
\mid \nu(x) \in \lambda(\A_{F'}^*F^*F^*_{\infty +})\},\]
where $\nu:\tilde{G}\to \G_m$ is the reduced norm. 
For any compact open subgroup $\tilde{K}$ of $
\tilde{G}(\A_{F,f})$, define a subgroup  $N(\tilde{K})$ of the group
of ideles $\A_{F'}^*$ as, 
\[ N(\tilde{K})=\{c\in \A_{F'}^*\mid \lambda(c)\in
F^*\nu(\tilde{K})\}.\]
By class field theory,  the subgroup $N(\tilde{K})$ defines an abelian
extension $F_{\tilde{K}}'$ of $F'$. Given $x\in \tilde{\mathcal G}_+$,
choose an element $c\in \A_{F'}^*$ such that $\lambda(c)/\nu(x)\in
F^*F^*_{\infty +}$. By the reciprocity morphism of class field,
\[ \mbox{rec}: \A_{F'}^*/F'^*\to \mbox{Gal}(F'^{ab}/F'),\]
we get an element $\sigma(x)\in  \mbox{Gal}(F'^{ab}/F')$ by the
prescription 
\[ \sigma(x)= \mbox{rec}(c^{-1}).\]

Associated to  $\tilde{K}$, we get a lattice
\[ \Gamma_{\tilde{K}}=\tilde{K}\tilde{G}_{\infty +}\cap
\tilde{G}{\Q}.\]
 Let  $K=\tilde{K}\cap G(\A_{F,f})$, and assume the following:\\

\noindent{\em Assumption (*).}\label{asscanmodel}
The natural
 inclusion of $\Gamma_K\subset \Gamma_{\tilde{K}}$ is an isomorphism
 modulo the centers. 

This assumption ensures that the quotient of  $M$ by the actions of
$\Gamma_K$ and $\Gamma_{\tilde{K}}$ are isomorphic. 

For sufficiently small $\tilde{K}$, the group $ \Gamma_{\tilde{K}}$
modulo the centre acts freely and discontinuously on the associated
symmetric space $M$. By the theory of canonical models, 
 the complex varieties $M_K$ acquire a canonical model
defined over the  abelian extension of $F_{\tilde{K}}'$, satisfying various
compatibility properties. We will refer to such varieties as `Shimura
varieties of quaternionic type'.  One of the principal properties that we
require is that there is an isomophism between the spaces,
\begin{equation}\label{conjcanon}
 M_K\simeq M_{K^x}^{\sigma(x)}.
\end{equation}

We  now give an `arithmetical archimedean analogue of Tate's
conjecture' given in \cite{PR}:
\begin{conjecture}\label{conj:aaat}
Suppose $X$ and $Y$ are Shimura varieties of quaternionic type of
dimension $s$,  which
are isospectral with respect to the deRham Laplacian acting on the
space of $p$-forms for $0\leq p <r$. Then the Hasse-Weil zeta
functions of $X$ and $Y$ are equal. 
\end{conjecture}

\begin{remark} A more general conjecture for curves is given in
  \cite{PR}, where we do not restrict ourselves to Shimura curves. But
  for such curves, it is not at all clear about the naturality of
  either the Riemannian metric nor of the arithmetic of such
  varieties. For Shimura curves and more generally for the class of
  spaces defined by quotients of symmetric spaces by congruent
  arithmetic lattices, both the arithmetic and the
  spectrum can be related to theory of automorphic forms, and hence it
  is reasonable to restrict our attention to such spaces.

For higher dimensional varieties, it is natural to add the spectrum of
deRham Laplacian acting on the space of smooth differential forms of
higher degrees, or even other natural elliptic self adjoint
differential operators. 
 But as we will see in Section \ref{sec:genconj}, it is more
natural to break up  this conjecture into two parts: one, relating the
the spectral side to   
representation equivalence of lattices; the second is to reformulate
the above conjecture on the spectral side in terms of representation
equivalence of lattices. 
\end{remark}

The first striking evidence for this conjecture that leads us to
believe in the above conjecture, comes from the
following theorem of A. Reid \cite{Re}:
\begin{theorem}[Reid]
Suppose $X$ and $Y$ are Shimura curves of quaternionic type,
associated respectively to quaternion division algebras $D_X, ~D_Y$
defined respectively over totally real number fields $F_X$ and
$F_Y$. Assume that $X$ and $Y$ are isospectral with respect to the
hyperbolic Laplacian acting on the space of smooth functions. Then
\[ F_X=F_Y \quad \mbox{and} \quad D_X=D_Y.\]
\end{theorem}

The proof of Reid's theorem starts with they equivalent hypothesis
that $X$ and $Y$ are length isopectral, a consequence of the Selberg
trace formula. The length spectrum is then related to the arithmetic
of the quaternion division algebras. This later step
 has been generalised by T. Chinburg, E. Hamilton,
D. Long and A. Reid \cite{CHER}, and more generally 
 by G. Prasad and A. S. Rapinchuk
\cite{PRap}. 

\section{Representation equivalence of lattices} \label{Gassmann-Sunada}
In this section, we relate the representation theoretic approach to
isospectral questions.  We recall the
Gassmann-Sunada method: A triple $(G, ~H_1, ~H_2)$ consisting of a
finite group $G$ and two non-conjugate subgroups $H_1, ~H_2$ of $G$ 
is said to form
a {\em Gassmann-Sunada system} if it satisfies one of the following
equivalent conditions:

\begin{enumerate}
\item The regular representations of $G$ on $\C[H_1\backslash G]$ and
  $\C[H_2\backslash G]$ are equivalent. 

\item Any $G$-conjugacy class intersects $H_1$ and $H_2$ in the same
  number of elements. 

\item For any finite dimensional representation $V$ of $G$, the
  dimension of the spaces of invariants with respect to the subgroups 
$H_1$ and $H_2$ are equal. 

\end{enumerate}

\begin{example}
There exists  examples of Gassmann-Sunada systems. 
 For example suppose $H_1$
and $H_2$ are two non-isomorphic 
finite groups such that for any natural number $d$, the
number of elements of order $d$ in $H_1$ and $H_2$ are equal (this is
a necessary condition for $H_1$ and $H_2$ to be part of a
Gassmann-Sunada system as above). In particular, 
the cardinalities of  $H_1$ and $H_2$ are equal, say $d$. Then $(S_d,
~H_1, ~ H_2)$ forms a Gassmann-Sunada system, where $S_d$ denotes the
symmetric group on $d$ symbols.     
\end{example}

Gassmann proved the following theorem (see \cite{PR}, \cite{Per}):
\begin{theorem}[Gassmann]\label{thm:gassmann}
Suppose $K$ is a finite Galois extension of
the rationals with Galois group isomorphic to $G$. Suppose  $G$ belongs
to a  Gassmann-Sunada system $(G, ~H_1, ~H_2)$. Then the field of invariants
$K^{H_1}$ and $K^{H_2}$ are non-conjugate number fields having the
same Dedekind zeta function.
\end{theorem}

Sunada's fundamental observation  (\cite{S}) is 
 that an entirely analogous statement can be
transplanted to the spectral side:
\begin{theorem}[Sunada] \label{thm:sunada}
Suppose $(G, ~H_1, ~H_2)$ is a
Gassmann-Sunada system and assume that 
 $G$ acts freely and isometrically on a
compact Riemannian manifold $M$. Then the quotient spaces of $M$ by
the action of the groups $H_1$ and 
 $H_2$ are isospectral with respect to the metric induced from $M$. 
\end{theorem}

The proofs of the theorems of Gassmann and Sunada as well as the
validity of Question \ref{qn:naiveaat}  for the examples of pairs of
isospectral compact Riemann surfaces constructed using the Sunada
method, rest on the functorial nature of the Frobenius isomorphism
 i.e., given a $G$-space $V$, and $H$ a subgroup of the
finite group  $G$, the Frobenius isomorphism,
\[ V^H\simeq (V\otimes \C[H\backslash G])^G\]
is natural (see \cite{PR}) with respect to $G$-equivariant maps $\phi:
V\to W$ of $G$-spaces. Sunada's theorem follows by taking
$V=C^{\infty}(M)$, and   $\phi$ to be the Laplacian. Since $G$ acts by
isometries, the Laplacian is $G$-equivariant, and restricts to give
the Laplacian of the space $V^{H_i}\simeq C^{\infty}(M/H_i), ~i=1,
~2$.  

Sunada's condition for isospectrality can be extended to infinite
groups, to the context of continuous Lie group actions and discrete
subgroups of such Lie groups. 
Let $G$ be a Lie group, and let $\Gamma$ be a cocompact lattice in
$G$. The existence of a cocompact lattice implies that $G$ is
unimodular. Let $R_{\Gamma}$ denote the right regular representation
of $G$ on the space $L^2(\Gamma\backslash G)$ of square integrable
functions with respect to the projection of the Haar measure on the
space $\Gamma\backslash G$:
\[(R_{\Gamma}(g)f)(x)=f(xg)\quad f\in  L^2(\Gamma\backslash G), ~g,
x\in G.\]
As a $G$-space, $L^2(\Gamma\backslash G)$ breaks up as a direct sum of
irreducible unitary representations of $G$, with each irreducible
representation occuring with finite multiplicity. 

\begin{definition} Let $G$ be a Lie group and $\Gamma_1$ and
  $\Gamma_2$ be two co-compact lattices in $G$. The lattices
$\Gamma_1$ and $\Gamma_2$ are said to be {\em representation
equivalent in $G$}  if the regular representations $R_{\Gamma_1}$ and
$R_{\Gamma_2}$ of $G$ are isomorphic.
\end{definition}

We have the following generalization of Sunada's criterion for
isospectrality proved by DeTurck and Gordon in \cite{DG}:
\begin{proposition}\label{sunada}
Let $G$ be a Lie group  acting on the left as isometries of a
Riemannian manifold $M$.  Suppose $\Gamma_1$ and $\Gamma_2$ are
discrete, co-compact  subgroups of $G$ acting freely and properly
discontinuously on $M$, such that the quotients $\Gamma_1\backslash M$
and $\Gamma_2\backslash M$ are compact Riemannian manifolds. 

If  the lattices  $\Gamma_1$ and $\Gamma_2$ are representation equivalent in
$G$, then  $\Gamma_1\backslash M$ and $\Gamma_2\backslash M$ are
isospectral for the Laplacian acting on the space of smooth functions.
\end{proposition} 

The conclusion in Proposition \ref{sunada} can be strengthened to
  imply  `strong isospectrality', i.e., the spectrums coincide for
  natural self-adjoint  differential operators besides the Laplacian
  acting on functions. In particular, this implies the isospectrality
  of the Laplacian acting on $p$-forms.

The natural hermitian vector bundles on spaces
  of the form $\Gamma\backslash G/K$ should be the class of
  automorphic vector bundles. On the universal cover $G/K$, the
  automorphic vector bundles are the $G$-invariant hermitian vector
  bundles that are associated to representations of $K$ on some finite
  dimensional unitary vector space $V$. The centre $Z(\frak g)$ of the
  universal enveloping algebra acts on the space of smooth sections
  $C^{\infty}(\Gamma\backslash G, V)$ of
  these vector bundles. Since the image of $Z(\frak g)$ is generated 
by invariant self adjoint elliptic differential operators, it follows
  that if $\Gamma$ is a uniform lattice in $G$, that  for 
any character $\chi$ of $Z(\frak g)$ the $\chi$-eigenspace, 
\[ \{\phi\in C^{\infty}(\Gamma\backslash G, V)\mid z.\phi=\chi(z)\phi,
  \quad z\in Z(\frak g)\},\]
has finite dimension  $m(\chi, \Gamma, V)$.

 Proposition \ref{sunada} can be refined  to say that the two spaces
 $\Gamma_1\backslash M$ and  $\Gamma_2\backslash M$ are {\em strongly
 compatibly isospectral}, in the sense that for any $V, ~\chi$ as
 above, we have an equality of multiplicities,
\[ m(\chi, \Gamma_1, V)=m(\chi, \Gamma_2, V).\]

\subsection{Converse to Gassmann-Sunada criterion}
We see that representation equivalent lattices give rise to  
spectrally indistinguishable spaces. 
It is natural to ask whether the converse holds. When $G=PSL_2(\R)$,
basic facts from the representation theory of $G$, allow us to show
that the converse holds.  
 Given a hyperbolic compact Riemann surface
$X$, let
\[ \rho_X :\Gamma_X\to PSL(2,\R),\]
be an embedding of the fundamental group $\Gamma_X$ arising from the
uniformization of $X$ by the upper half plane. This map is independent
of the various choices made upto conjugation by an element of
$PSL(2,\R)$. A folklore observation is the following \cite{Pes2}:
\begin{theorem}\label{thm:convsunadpsl2}
Suppose $X$ and $Y$ are compact Riemann surfaces of genus at least
two.  Then $X$ and $Y$ are isospectral with respect to the hyperbolic
Laplacian acting on the space of smooth functions if and only if
the lattices $\rho_X(\Gamma_X)$and $\rho_Y(\Gamma_Y) $ are
representation equivalent in $PSL(2,\R)$: 
\[L^2(PSL(2,\R)/\Gamma_X)\simeq L^2(PSL(2,\R)/\Gamma_Y),\]
as $PSL(2,\R)$-modules.
\end{theorem}
The proof that the spectrum determines the
representation type of $\Gamma_X$ in $PSL(2,\R)$ 
follows from the classification of the
irreducible unitary representations of $PSL(2,\R)$. 

\begin{remark} Maass showed that a theory of Hecke operators can
  be defined on the spaces of  eigenfunctions of the
  Laplacians invariant by the lattice, paving the way for the work of Selberg
  and Langlands. 
\end{remark}

\begin{remark}
One of the main reasons for considering the hyperbolic Laplacian for
Riemann surfaces 
is that the spectrum of the hyperbolic Laplacian can be related
to the  representation theory of the isometry group $PSL(2,\R)$ of the
universal covering space.If we consider other natural metrics on a
compact Riemann surface, then we can no longer expect any relationship
with representation theory. 
\end{remark}

For groups apart from $PSL_2(\R)$, the converse direction from the
spectrum to the representation theory is not well understood. 
For compact quotients of hypebolic spaces, H. Pesce \cite{Pes1} has
shown that if two such spaces are  strongly isospectral, then 
the corresponding lattices are
representation equivalent in the group of isometries of the hyperbolic
space. For more general groups, we formulate
later a possible approach to   a converse of  the generalized Sunada
criterion in the context of compact locally symmetric spaces (see Conjecture 
\ref{conj:convsunad}).

\section{Adelic conjugation of lattices}
In this section, we provide some evidence for the conjectures raised
in the previous sections. In the process, we construct new examples of
isospectral but non-isometric spaces.

Let $F$ be a number field and let $D$ be  an indefinite
 quaternion division algebra 
over $F$. We assume that there is at least one  finite place $v_0$
of $F$, at which $D$ is ramified. Let  $K$ be   a compact open
subgroup of $G(\A_{F,f})$ having a factorisation of the form, 
\[ K=K_0K',\]
where $K_0$ is a compact,  open subgroup of $G(F_{v_0})$,  and is an invariant 
subgroup of $D^*_{v_0}$. The group $K'$ has no $v_0$ component, i.e.,
for any element $x\in K'$,
the $v_0$ component $x_{v_0}=1$.  Let 
\[ \Gamma_K=G_{\infty}K\cap G(F),  \]
be the co-compact lattice in $G_{\infty}$, and let
$M_K=\Gamma_K\backslash M$ be the quotient space of $M$ by
$\Gamma_K$. 
 
The following theorem is proved in \cite{Ra}, and can be considered
as a geometric version of the concepts of $L$-indistinguishablity due
to Labesse, Langlands and Shelstad: 
\begin{theorem} \label{adelicconjtheorem}
With the above notations and assumptions, for any element 
$x\in GL_1(D)(\A_{F,f})$, the lattices $\Gamma_K$  and
$\Gamma_{K^x}$ are representation equivalent.  
\end{theorem}

The following corollary gives examples of isospectral but
non-isometric compact Riemannian manifold, generalizing the class of
examples constructed by Vigneras \cite{V}, 
reflecting essentially the failure of strong
approximation in the adjoint group:
\begin{corollary}\label{isospnotisom} With the notation as in Theorem
  \ref{adelicconjtheorem}, assume further that $K$ is small enough so
  that $\Gamma_K$ is torsion-free. 
 Let $\widehat{N(K)}$ denote the normalizer of the image of  $K$ in
$PGL_1(D)(\A_{F,f})$. Choose an element $x\in \tilde{G}(\A_f)$ such that
it's projection to $PGL_1(D)(\A_{F,f})$ does not belong to the set 
$\widehat{N({K})}PGL_1(D)(F)$ (such
elements exist by the failure of strong approximation in
$PGL_1(D)$). 

Then $X_K$ and
$X_{K^x}$ are  strongly isospectral, but are not isometric. 
\end{corollary} 

\begin{remark}
Vigneras works with the length spectrum rather than the representation
theoretic context in which the above theorem is placed. The relation
between the Laplacian spectrum and the length spectrum holds only for
hyperbolic surfaces and three folds. Further, she has to compute the
length spectrum and show that indeed the lattices are length
isospectral. For this reason, she has to restrict attention to $K$
which come from maximal orders and use theorems of Eichler computing
the number of elements having a given trace. 
\end{remark}

We also obtain the following corollary
providing more evidence in support of the Conjecture \ref{conj:aaat}:
\begin{corollary} \label{cor:naiveatfalse}
 Let  $F$ be  a totally real number field and 
$\tilde{K}$ be a compact open subgroup of $ \tilde{G}(\A_{F,f})$. Assume
  that $\tilde{K}$ satisfies the hypothesis of Theorem
  \ref{adelicconjtheorem} and is such  that the lattice
  $\Gamma_{K}$ is torsion-free and satisfies Assumption (*). 

 Then the spaces $M_{K^x}$ for $x\in
  \tilde{G}(\A_{F,f})$ are isospectral and have the same Hasse-Weil zeta
  function for the canonical model defined by Shimura.

If $D$ is ramified at all real places except one, then the Jacobians
of $M_{K}$ and $M_{K^x}$ are not isogenous but are 
conjugate by an automorphism of $\bar{\Q}$.
\end{corollary} 

This gives us the counterexamples to Question \ref{qn:naiveaat}.
It is tempting to conjecture that the last statement of the foregoing
corollary will be the exception to Question \ref{qn:naiveaat}:  if two
  compact Riemann surfaces are isospectral, then the Jacobian of one
  is isogenous to a conjugate of the Jacobian of the other by an
  automorphism $\sigma\in {\rm Aut}(\C/\Q)$, where $\sigma$ preserves
  the spectrum of the Riemann surface. If morever $\sigma$ is not
  identity or of order two, then the pair of Riemann surfaces  arise
  from an arithmetical context, i.e., are Shimura curves as considered
  in this paper.

\section{L-equivalence}
The proof of Theorem \ref{adelicconjtheorem} given in \cite{Ra}, 
uses the theory of
$L$-indistinguishable automorphic representations initiated by  Langlands,
in particular the computation of the mutltiplicity  
by Labesse and Langlands with which a
representation of $SL(1,D)(\A_F)$ occurs in the space of automorphic
representations of $SL(1,D)(\A_F)$ \cite{Cor}, \cite{LL}.
An examination of the proof of Theorem \ref{adelicconjtheorem}
indicates that it yields a more general formulation, which will form
the basis of reformulating the above conjectures in a more general
framework. 

From the viewpoint of Langlands theory, it is more natural to classify
the $L$-packets of irreducible  representations (\cite{A, Cor, KZ, L,
Sh}).   We will assume
that there exists a suitable notion of $L$-packets.
We do not attempt  out here the definition of either the Langlands or
Arthur $L$-packets, nor what is the correct notion of $L$-packets that
we require for the relationship between the spectrum and arithmetic to
hold. A coarse expectation
that two representations are $L$-equivalent (or $L$-indistinguishable
or equivalently belong to the same $L$-packet) is that the
$L$-functions and $\epsilon$-factors  attached to them are equal. 
It is expected that
$L$-packets have finite cardinality. For representations of real
Lie groups, a property that  members of the same $L$-packet share is
that they all have the same infinitesimal character, i.e., they cannot
be distinguished by means of spectral data alone.

\begin{example} If $F$ is a local field, the $L$-packets of
  representations of  $SL(2,F)$ (or any of it's inner forms
  $SL(1,D)(F)$)
consists of those representations that 
are conjugate to each other by an element of $GL(2,F)$ (resp. by
$GL(1,D)(F)$). Thus the
discrete series $\{\pi_k^+. ~\pi_k^-, ~k\geq 2\}$ 
constitute a single $L$-packet
for $SL(2,\R)$. 
\end{example}

Let $F$ be a global field. 
Two representations of $SL(2,\A_F)$ (resp. of $SL(1,D)(\A_F)$)
are said to be {\em $L$-equivalent} if
their local components at each place of $F$ are $L$-indistinguishable
as representations of the local group. Equivalently, they are
conjugate by an element of $GL(2,\A_F)$ (resp. $GL(1,D)(\A_F)$).

 Given an irreducible unitary representation $\pi$ of a semisimple Lie
group $G$ (real or $p$-adic or adelic), 
we denote by $[\pi]$, the $L$-packet of $G$ containing
$\pi$. For an automorphic
$\pi$ of $G(\A)$, we denote by $\pi^K$ the
space of $K$-fixed vectors of the space underlying the representation
$\pi$ and by $m_0(\pi)$ the (finite) multiplicity with which
$\pi$ occurs in the space of cusp forms. We 
define a notion of automorphic equivalence of lattices:

\begin{definition} Let $G$ be a reductive algebraic group over $F$ and
  $K$ be a compact open subgroup of $G(\A_{F,f})$. 
The {\em cupidal automorphic spectrum} associated to $(G(\A_F), K)$, 
is the collection of
cuspidal automorphic representations $\pi$ of $G(\A_F)$ counted with
multiplicity $m(\pi, K)$, where 
\[ m(\pi, K)= m_0(\pi)\mbox{dim}(\pi^K).\]

The {\em $L$-cuspidal spectrum} associated to  $(G(\A_F),
  K)$ is the collection of cuspidal $L$-packets of 
$[\pi]$ of $G(\A_F)$
counted with multiplicity $m([\pi], K)$ defined as, 
\[ m([\pi], K)=\sum_{\eta\in [\pi]}m(\eta,K).\]

\end{definition}

\begin{definition} Let $G_1$ and $G_2$ be semisimple groups over $F$. 
Suppose $K_1$ and $K_2$  are compact open subgroups of $G_1(\A_{F,f})$
and  $G_2(\A_{F,f})$ respectively. Then the pair  $(G_1(\A_{F,f}),K_1)$
and  $(G_2(\A_F), K_2)$ are said to be cuspidally equivalent (resp. 
$L$-cuspidally equivalent) if there exists an isomorphism 
\[ \phi:~G_1(\A_{F})\to G_2(\A_{F})\]
such that for any cuspidal
representation  $\pi$ of $ G_2(\A_{F})$, the
multiplicities $m(\pi,K_2)$ and $m(\pi\circ\phi, K_1)$
(resp. $m([\pi],K_2)$ and $m([\pi\circ\phi], K_1)$) are equal. 
\end{definition}
 
The reasons for the above definitions stem from the following
strengthening of Theorem \ref{adelicconjtheorem}:
\begin{theorem}
With the same hypothesis as in Theorem \ref{adelicconjtheorem},  the
lattices $K$ and $K^x$ are $L$-equivalent in $SL_1(D)(\A_F)$. 
\end{theorem}
The proof of  Theorem
\ref{adelicconjtheorem} given in \cite{Ra} carries over to prove this
more general theorem.  

\section{From Arithmetic to Spectrum} \label{arithtospectrum}
Suppose $X$ is a projective  Shimura curve. We take  for a tentative
definition of  the arithmetic of $X$ the data given by 
Hasse-Weil zeta functions attached to `natural' local systems on $X$,
i.e., those local systems which arise by restriction to the lattice 
from algebraic representations of the group $\tilde{G}=GL_1(D)$. The
irreducible representations are
essentially given by the symmetric powers of the standard
representation of $GL(2)$.  We are then
led to ask the converse to the above questions that the spectrum
determines the arithmetic:
\begin{question}
Suppose $X$ and $Y$ are Shimura curves, such that the Hasse-Weil 
zeta functions 
associated to the local systems coming from the $n$-th symmetric power
of the standard representation of $GL(2)$ are equal for all $n\geq
0$. 

Are $X$ and $Y$ isospectral with respect to the hyperbolic
Laplacian acting on the space of functions? 
\end{question}

There is one context however where the arithmetic associated to a pair
of varieties is expected to be the same. 
 For this, we assume that $X$ is 
associated to congruence lattices $K_X\subset
G(\A_{F,f})$ as described in Section \ref{section:aaat}.
 In this case, the canonical model of
the varieties are defined over an abelian extension of the reflex
field  $F'$. Now it
follows from the theory of canonical models for algebraic local
systems as above, that if $\sigma\in \mbox{Gal}(F'^{ab}/F')$, then 
the conjugate $Y=X^{\sigma}$ is first of all again a Shimura curve
associated to $D$; secondly, it has the same arithmetic as $X$ in the
sense described above. 

In
this case, we have the following theorem that Galois conjugation
preserves the spectrum:
\begin{theorem}
With notation as in Section \ref{section:aaat}, assume
  that $\tilde{K}$ satisfies the hypothesis of Theorem
  \ref{adelicconjtheorem} and is such  that the lattices
  $\Gamma_{K^x}$ are torsion-free and satisfies  Assumption (*). 

For any  $\sigma\in \mbox{Gal}(F'^{ab}/F')$, the spaces 
$X$ and $X^{\sigma}$ are isospectral for the hyperbolic Laplacian
acting on functions.
\end{theorem}

The proof of this theorem is based on the observation that by the
properties of canonical models, any $\sigma$ as above can be written
as $\sigma(x)$ for some $x\in GL_1(D)(\A_{F,f})$, and then appealing
to Theorem \ref{adelicconjtheorem}. In some sense, the `correct
interpretation'  of   Theorem \ref{adelicconjtheorem} is given by the
above theorem:  {\em Galois conjugation over the reflex field 
preserves the spectrum}. 

\section{General conjectures} \label{sec:genconj}
We try to generalize the discussion so far to the general case of
quotients of symmetric spaces defined by congruence arithmetic
lattices. The relationship between the spectrum and arithmetic can be
broken into two steps: one, the relationship between the spectrum and
$L$-equivalence of uniform lattices in real Lie groups. The other
aspect is to relate $L$-equivalence of lattices in the real points of
a reductive or semisimple algebraic group to suitable notions of
arithmetic, i.e., to relate to more global aspects. 

\subsection{Spectrum to arithmetic}
 Suppose $\Gamma$ is a cocompact
lattice in a real Lie group $G$. We have a direct sum decomposition, 
\[  L^2(\Gamma\backslash G)=\oplus_{\pi\in \hat{G}}m(\pi,
\Gamma)\pi,\]
where $\hat{G}$ denotes the set of equivalence classes of irreducible
unitary representations of $G$, and $m(\pi, \Gamma)$ the (finite)
multiplicity with which $\pi$ occurs in $ L^2(\Gamma\backslash
G)$. Define the multiplicity $m([\pi],\Gamma)$ of an $L$-packet
$[\pi]$ by,
\[m([\pi],\Gamma):=\sum_{\eta\in [\pi]} m(\eta,\Gamma).\]
Define two uniform lattices $\Gamma_1,~\Gamma_2$ to be $L$-equivalent
or $L$-indistinguishable in $G$, if for any $\pi\in \hat{G}$, 
\[  m([\pi],\Gamma_1)=m([\pi],\Gamma_2).\]
With this definition, we can expect the following conjecture to hold:
\begin{conjecture}\label{conj:convsunad}
Let $G$ be a semisimple Lie group, $K$ a maximal compact subgroup of
$G$ and  $\Gamma_1,~\Gamma_2$ two torsion-free uniform lattices in
$G$. 

Then  $\Gamma_1$ and $\Gamma_2$ are $L$-equivalent in $G$ if and
only if the associated compact locally symmetric spaces
$\Gamma_1\backslash G/K$ and $\Gamma_2\backslash G/K$ are compatibly
strongly isospectral as defined in Section \ref{Gassmann-Sunada}.
\end{conjecture}

One can replace the isospectrality on forms by  weaker
assumptions: for instance, isospectrality 
for the deRham Laplacian on $p$-forms for all $0\leq p\leq  \mbox{dim}(G/K)$.
Or for the compact quotients
of symmetric spaces of noncompact type, it is even tempting to ask
whether isospectrality on functions is sufficient to guarantee the
rest of the implication. 

\begin{remark} When the lattices are no longer cocompact, one can
  define a similar notion of two lattices being {\em cuspidally
    L-equivalent} by working with the cuspidal spectrum of
  $L^2(\Gamma\backslash G)$. 
\end{remark}

The relationship between $L$-equivalence and
representation equivalence of lattices is not clear.  
For instance, the following
question can be raised:
\begin{question} Does there exist $L$-equivalent lattices in
 which are not representation equivalent?
\end{question}

\begin{remark}
For $SL(2,\R)$, it seems that Hodge duality will ensure that 
 $L$-equivalent lattices are representation equivalent. 
\end{remark}

\begin{remark} To start the discussion rolling from the spectral side,
  a first question to ask is whether any space isospectral to a
  compact locally symmetric space is itself locally symmetric? If we
  assume strong isospectrality, then this is proved by Gilkey
  \cite{Gi}.  A suitable assumption will be to impose isospectrality on
  forms, but for the compact quotients of non-compact symmetric spaces
  it would be interesting to know whether just the spectrum on
  functions will suffice.

Assuming the truth of the above question, the natural sequel to it is
to say that the spectrum determines the symmetric space, in other
words the isometry group of the universal cover. 

Such results can be considered as spectral analogues of Langlands
conjectures on conjugation of Shimura varieties, that the Galois
conjugate of a Shimura variety is again a Shimura variety. 
\end{remark}

We now generalize the conjecture that the spectrum determines the
arithmetic:
\begin{conjecture}\label{conj:gen:spectrumtoarith}
Let $G_1, ~G_2$ be semisimple algebraic groups over number field $F_1$
and $F_2$ respectively. Let $K_i\subset G_i(\A_{F,f}), ~i=1,2$ be compact open
subgroups, and $\Gamma_{K_i}$ be the corresponding arithmetic lattices in
$G_{i, \infty}$. Suppose $G_{1, \infty}\simeq
G_{2, \infty}=G_0$, and that  $\Gamma_{K_1}$
and $\Gamma_{K_2}$ are $L$-cuspidal, equivalent lattices in $G_0$. 

Then  the adele groups $G_1(\A)$ and $G_2(\A)$
are isomorphic, and the $L$-cuspidal spectrums of
$(G_1(\A_F), K_1)$ and $(G_2(\A_F), K_2)$ are equal. 
\end{conjecture}

\begin{remark} The results of Prasad-Rapinchuk \cite{PRap} show that if
the spaces corresponding to the two lattices satisfy a `weak
commensurability' property, then the
corresponding adele groups are isomorphic.
\end{remark}
 
The following $GL_1$-aspect of the above conjecture can
  be made: 
\begin{conjecture} Suppose $F_1, ~F_2$ are two totally real number
  fields of degree $r$ over $\Q$ are 
such that the unit groups $U_{F_1}$ and $U_{F_2}$ are spectrally
commensurable as lattices in the natural embedding into
$\R^{r-1}$. Then are $F_1$ and $F_1$ arithmetically equivalent, i.e.,
do they have the same Dedekind zeta function? 
\end{conjecture}

It is known (see \cite{Per}) that if two number fields are
arithmetically equivalent then the unit groups have the same rank. The
methods of \cite{PR} will also directly show that the lattices are
commensurable with respect to the natural embeddings, since it is
known that arithmetically equivalent number fields arise from the
Gassmann construction. But for the converse direction, we need to
restrict our attention to totally real number fields. It is for this
reason that we have restricted our attention mostly to semisimple
groups, although it is quite natural to frame the conjectures in the
wider context of reductive algebraic groups.

\subsection{From arithmetic to spectrum}
We now consider the possible definitions of arithmetic not only 
in the context of higher dimensional Shimura varieties but also for
more general locally symmetric spaces arising from congruent arithmetic
lattices.

Let  $G$ be a semisimple algebraic group over a number
field $F$, and $K$ be a compact open subgroup of  
$G(\A_{F,f})$. One possible definition for the arithmetic associated
to $(G(\A_F),K)$ is to consider the cuspidal $L$-spectrum 
attached to it. A natural question is the following:
\begin{conjecture}\label{conj:gen:autarithtospectrum}
Let $G_1, ~G_2$ be semisimple algebraic groups over number field $F_1$
and $F_2$ respectively. Let $K_i\subset G_i(\A_{F,f}), ~i=1,2$ be compact open
subgroups, and $\Gamma_{K_i}$ be the corresponding arithmetic lattices in
$G_{i, \infty}$. Suppose that the adele groups $G_1(\A)$ and $G_2(\A)$
are isomorphic, and that the  $L$-cuspidal spectrums of
$(G_1(\A_F), K_1)$ and $(G_2(\A_F), K_2)$ are equal. 

Then $\Gamma_{K_1}$
and $\Gamma_{K_2}$ are $L$-equivalent lattices in $G_{1, \infty}\simeq
G_{2, \infty}$.
\end{conjecture}
 
This conjecture seems obvious when $F_1=F_2$, and 
$G_1$ and $G_2$ are isomorphic. However Lubotzky, Vishne and Samuels
\cite{LSV} have constructed non-trivial examples satisfying the
hypothesis of the above conjecture, concluding thereby that there
exists non-commensurable arithmetic lattices in $PGL_d(\R)$ which are
representation equivalent. The above conjecture should yield to
suitable generalizations of Jacquet-Langlands type comparison between
automorphic representations on inner forms of algebraic groups. 

We now try to  define a suitable notion of arithmetic for
general groups not giving raise to hermitian symmetric spaces, 
which closely reflects the customary meaning of arithmetic of Shimura
varieties as discussed in Section \ref{arithtospectrum}.

We recall that an automorphic representation $\pi$ of $G(\A_F)$
is cohomological if it's infinity type $\pi_{\infty}$ is a
cohomological representation of $G(F\otimes_{\Q} \R)$. For example,
the discrete series representations $\pi_k^{\pm}, ~k\geq 2$ are
representations of cohomological type for $SL(2,\R)$ (together with
the trivial representation and the Steinberg representation they
constitute all the representations of $SL(2,\R)$ with cohomology). 

It follows from the congruence relation of Eichler-Shimura and the
work of Ihara, Langlands and others (\cite{PS}, \cite{Cor}), 
 that in case of Shimura varieties
the automorphic representations that describe the Hasse-Weil zeta
functions of  natural local systems on such spaces, are cohomological.   
\begin{definition}
The {\em cohomological $L$-cuspidal  spectrum} 
associated to  $(G(\A_F),
  K)$ is the collection of $L$-packets of cohomological 
automorphic representations of $G(\A_F)$
counted with multiplicity $m([\pi], K)$.
\end{definition}

It's this definition, that we would like to consider as the
suitable candidate for the notion of arithmetic of the space
$\Gamma_K\backslash G_{\infty}/K_{\infty}$. However there are a couple
of problems with this definition: one, is that computations done for
general groups seem to indicate a sparsity of cohomological
representations; the other is that we will have to consider
non-unitary cohomological representations also in the above
definition.

We now generalize the conjecture that the arithmetic of Shimura curves
determines the spectrum to arbitrary locally symmetric
spaces:
\begin{conjecture}\label{conj:gen:arithtospectrum}
Let $G_1, ~G_2$ be semisimple algebraic groups defined over number fields $F_1$
and $F_2$ respectively. Let $K_i\subset G_i(\A_{F,f}), ~i=1,2$ be compact open
subgroups, and $\Gamma_{K_i}$ be the corresponding arithmetic lattices in
$G_{i, \infty}$. Suppose that the adele groups $G_1(\A)$ and $G_2(\A)$
are isomorphic, and that cohomological $L$-cuspidal spectrums of
$(G_1(\A_F), K_1)$ and $(G_2(\A_F), K_2)$ are equal. 

Then $\Gamma_{K_1}$
and $\Gamma_{K_2}$ are $L$-equivalent lattices in $G_{1, \infty}\simeq
G_{2, \infty}$.
\end{conjecture}

A particular consequence of Conjectures \ref{conj:gen:spectrumtoarith}
and \ref{conj:gen:arithtospectrum} is
that the cuspidal cohomological  $L$-spectrum should determine the full
cuspidal $L$-spectrum.

\begin{remark} If we consider the example of $SL(1,D)$, where $D$ is
  an indefinite quaternion division algebra over a totally
  real number field, the arithmetic
  is concentrated on those automorphic representations whose
  archimedean component is a discrete series. The spectrum on the
  other hand gives information about those automorphic representations
  with archimedean component an  
unramified principal series
  representation. 
\end{remark}

\end{document}